\documentclass[a4paper,11pt]{article}
\usepackage{amsmath,amscd,amsfonts,amssymb,epsf,latexsym}

\makeatletter

\long\def\@makefnt#1{\parindent 1em\noindent
            \hb@xt@1.8em{\hss\@textsuperscript{}}#1}
\long\def\@ftntext#1{\insert\footins{%
    \reset@font\footnotesize
    \interlinepenalty\interfootnotelinepenalty
    \splittopskip\footnotesep
    \splitmaxdepth \dp\strutbox \floatingpenalty \@MM
    \hsize\columnwidth \@parboxrestore
    \color@begingroup
      \@makefnt{%
        \rule\z@\footnotesep\ignorespaces#1\@finalstrut\strutbox}%
    \color@endgroup}}%
\def\subjclass#1{%
  \@ftntext{2000 {\itshape Mathematics Subject Classification.}\enspace #1.}}
\def\keywords#1{%
  \@ftntext{{\itshape Key words and phrases.}\enspace #1.}}
\makeatother

\def\moins{\raise 1pt\hbox{{$\scriptstyle -$}}}
\def\plus{\raise 1pt\hbox{{$\scriptstyle +$}} }

\newtheorem{theorem}{Theorem}

\newtheorem{remark}[theorem]{Remark}

\newtheorem{example}[theorem]{Example}
\newtheorem{note}[theorem]{Note}

\def\proof{\noindent{\bf Proof.\ }}

\def\qed{~\hbox{$\Box$}}

\def\Aut{\mathop{\rm Aut}}
\def\cJ{\mathop{\rm J}}

\def\codim{\mathop{\rm codim}}

\def\rank{\mathop{\rm rank}}
\def\dim{\mathop{\rm dim}}
\def\Sym{\mathop{\rm Sym}}

\def\cT{{\mathcal T}}
\def\cJ{{\mathcal J}}
\def\cS{{\mathcal S}}

\begin{document}

\title{\bf Positivity of Schur function expansions\\
of Thom polynomials}

\author{Piotr Pragacz\thanks{Research supported by a KBN grant and
by the Humboldt Stiftung during the stay at the MPIM in Bonn.}\\
\small Institute of Mathematics of Polish Academy of Sciences\\
\small \'Sniadeckich 8, 00-956 Warszawa, Poland\\
\small P.Pragacz@impan.gov.pl
\and
Andrzej Weber\thanks{Research supported by the KBN grant 1 P03A 005 26.}\\
\small Department of Mathematics of Warsaw University\\
\small Banacha 2, 02-097 Warszawa, Poland\\
\small aweber@mimuw.edu.pl}

\date{(10.05.2006; revised 12.11.2006)}



\maketitle

\begin{abstract}
Combining the approach to Thom polynomials via classifying
spaces of singularities with the Fulton-Lazarsfeld theory of 
cone classes and positive polynomials for ample vector bundles,
we show that the coefficients of the Schur function expansions of
the Thom polynomials of stable singularities are nonnegative with
positive sum.
\end{abstract}

\section{Introduction}
The global behavior of singularities\footnote{In the present paper we
study complex singularities.} is governed by their {\it Thom
polynomials} (cf. \cite{T}, \cite{AVGL}, \cite{Ka}, \cite{Rim}).
As these polynomials are quite complicated even for ``simplest''
singularities, it is important to study their structure.
There is a recent attempt to present Thom polynomials via their
{\it Schur function expansions} (cf. \cite{FK}, \cite{P}, \cite{P1}) 
instead of using the ``traditional'' basis of monomials in Chern classes.

In the present paper, we study the Schur function expansions of Thom
polynomials from a ``qualitive'' point of view.
Contrary to \cite{Rim}, \cite{P}, \cite{P1}, where the Sz\"ucs-Rimanyi 
approach via symmetries of singularities was used, we follow here the Kazarian
approach \cite{Ka} to Thom polynomials.
In fact, both approaches rely on suitable ``classifying spaces of singularities''.
We substitute the jet automorphism group by the group of the linear
transformations $GL_m\times GL_n$. This allows one to extend the
definition of Thom polynomials for the maps $f: M\to N$ of complex manifolds
to pairs of vector bundles. It is convenient to pass to  homotopy theory,
where each pair of bundles can be pulled back from the universal pair
of bundles on $BGL_m\times BGL_n$.

We apply the Fulton-Lazarsfeld theory of cone classes and
positive polynomials for ample vector bundles \cite{FL}, and deduce
{\it nonnegativity} of the coefficients in the Schur function expansions
of the Thom polynomials of the singularities stable under suspension,
with positive sum.

This ``positivity'' was previously checked for a number of singularities:
by Thom \cite{T} for $A_1(r)$,
by Feher and Komuves \cite{FK} for some second order Thom-Boardman
singularities,
by the first author \cite{P}, \cite{P1} for $I_{2,2}(r)$, $A_3(r)$, 
and for $A_i(r)$
under the aditional assumption that $\Sigma^j(f)=\emptyset$ for $j\ge 2$,
by the first author and Ozturk for Thom polynomials from \cite{Rim},
by the second author for the Thom polynomials (from \cite{Ka})
of singularities of functions,
and by Ozturk \cite{O} for $A_4(3)$, $A_4(4)$. \footnote{We use here the
notation from \cite{P}, \cite{P1}. The calculations
in the last three cases used extensively ACE \cite{V}.} Some of these
examples are listed in the last section.

The paper is a revised version of the MPIM Bonn Preprint 2006-60.

\section{Thom polynomials}\label{Thom}
Fix $m,n,k \in {\bf N}$. We denote by $\Aut_n$ the group of
$k$-jets of automorphisms of $({\bf C}^n,0)$, and by
$\cJ=\cJ(m,n)$ the space of $k$-jets of functions $({\bf
C}^m,0)\to ({\bf C}^n,0)$ \footnote{Though these objects depend on
$k$, we omit ``$k$'' in the notation. This will happen also to
other objects introduced later.}.

Moreover, we set
$$
G:={\Aut}_m\times {\Aut}_n\,.
$$
Consider the classifying principal $G$-bundle $EG\to BG$, i.e.
a contractible space $EG$ with a free action of the group $G$, and
define
$$
\widetilde{\cJ}=\widetilde{\cJ}(m,n)=EG\times_G \cJ\,.
$$
Let $\Sigma \subset \cJ$ be an analytic closed $G$-invariant
subset, which we shall call a ``class of singularities''.
For a given class of singularities $\Sigma$, set
$$
\widetilde{\Sigma}=EG\times_G \Sigma\subset\widetilde{\cJ}
$$
and denote by ${\cT}^{\Sigma} \in
H^{2\codim(\Sigma)}(\widetilde{\cJ}, {\bf Z})$ the dual
class of $[\widetilde{\Sigma}]$ \footnote{One may approximate 
$EG\to BG$ by a sequence of $G$-bundles over
finite dimensional manifolds $(EG)_N\to (BG)_N$, where $N\to \infty$. Then
``${\cT}^{\Sigma}$ is the dual class of $[\widetilde{\Sigma}]$''
means that for any $N$ the image of ${\cT}^{\Sigma}$ in
$H^{\bullet}((BG)_N, {\bf Z})\cong H^{\bullet}((EG)_N\times_G\cJ,
{\bf Z})$ is Poincar\'e dual to the class $[(EG)_N\times_G\Sigma]$
in $H_{\bullet}^{BM}((EG)_N\times_G\cJ, {\bf Z})$.}.
Since
$$
H^{\bullet}(\widetilde{\cJ}, {\bf Z})
\simeq H^{\bullet}(BG, {\bf Z})\simeq H^{\bullet}(BGL_m\times BGL_n, {\bf Z})\,,
$$
${\cT}^{\Sigma}$ is identified with a polynomial in $c_1,\ldots, c_m$
and $c_1',\ldots, c_n'$ which are the Chern classes of universal bundles
$R_m$ and $R_n$ on $BGL_m$ and $BGL_n$:
$$
{\cT}^{\Sigma}={\cT}^{\Sigma}(c_1,\ldots,c_m,c'_1,\ldots,c'_n)
$$
(here, and in the following, we omit the pull back in the notation).
This is a classical {\it Thom polynomial}.
Speaking slightly informally, given a general map $f: M \to N$ of
smooth varieties of corresponding dimensions $m$ and $n$, the Thom polynomial
$$
{\cT}^{\Sigma}(c_1(M),\ldots, c_m(M),c_1(N),\ldots, c_n(N))
$$
evaluates the dual class of the set where $f$ has
singularity ``of the class $\Sigma$''. A precise version of this
statement is a content of the Thom theorem \cite{T} (see also \cite{Ka},
Theorem 1 and \cite{Rim}, Sect.~6).

The {\it suspension}
$$
\cS:\cJ(m,n)\hookrightarrow \cJ(m+1,n+1)
$$
allows one to increase the dimension of the source and the target
simultaneously: with the local coordinates $x_1,x_2,\ldots$ for the source
and a function $f=f(x_1,\ldots,x_m)$, the jet
$({\cS}f)\in\cJ(m+1,n+1)$ is defined by
$$
({\cS}f)(x_1,\dots,x_m,x_{m+1}):=(f(x_1,\dots,x_m),x_{m+1})\,.
$$
Suppose that the class of singularities $\Sigma$ is
{\it stable under suspension}. By this we mean that it is a member
$\Sigma_0=\Sigma$ of a family
$$
\{\Sigma_r\subset\cJ(m+r,n+r)\}_{r\ge 0}
$$
such that
$$
\Sigma_{r+1}\cap\cJ(m+r,n+r)=\Sigma_r
$$
and
\begin{equation}\label{sp}
{{\cT}^{\Sigma_{r+1}}}_{|H^{\bullet}(BGL_{m+r}\times BGL_{n+r}, {\bf Z})}
={\cT}^{\Sigma_r}\,.
\end{equation}
This means that if we specialize
$$
c_{m+r+1}=c'_{n+r+1}=0
$$
in the polynomial ${\cT}^{\Sigma_{r+1}}$, we obtain the polynomial
${\cT}^{\Sigma_r}$. If the class $\Sigma$ is stable under
V-equivalence (cf. \cite[\S I.6.5]{AGV}) then it is stable in our sense.

The theorem of Thom has the following refinement due to Damon
\cite{D} for a class of singularities $\Sigma$ which is stable
under suspension: ${\cT}^{\Sigma}$ is {\it supersymmetric}, i.e.
is a polynomial in
$$
(1+c_1+\ldots+c_m)/(1+c'_1+\ldots+c'_n) \ \ \ \hbox{where} \ \ \
i=1,2,\ldots\,.
$$
In other words, for a general map $f: M\to N$,
$$
{\cT}^{\Sigma}(c_1(M),\ldots, c_m(M),c_1(N),\ldots, c_n(N))
$$
is a polynomial in
$$
c_i(TM\moins f^*TN)=[c(TM)/c(f^*TN)]_i \ \ \ \hbox{where} \ \ \
i=1,2,\ldots\,
$$
(here, $TM\moins f^*TN$ is a virtual bundle).
Cf. also \cite[Theorem 2]{Ka}.

\section{Schur functions expansions}

Given a partition $I=(i_1,i_2,\ldots,i_l)\in {\bf N}^l$, 
where $0\le i_1\le i_2\le \ldots\le i_l$,
and vector bundles $E$ and $F$ on some variety $X$,
the {\it Schur function}\footnote{Usually this family of functions is called
``super Schur functions'' or ``Schur functions in difference of bundles'';
the classical Schur functions ${\cal S}_I$ will be used in 
Theorem \ref{th2} in the next section.}
$S_I(E - F)$ is defined by the following determinant:
\begin{equation}\label{schur}
S_I(E - F):= \Bigl|
     S_{i_p+p-q}(E - F) \Bigr|_{1\leq p,q \leq l}\,,
\end{equation}
where the entries are defined by the expression
\begin{equation}\label{seg}
\sum S_i(E - F) =\prod_b(1 - b)/\prod_a(1 - a)\,.
\end{equation}
Here, the $a$'s and $b$'s are the Chern roots of $E$ and $F$
and the LHS of Eq. (\ref{seg}) is the {\it Segre class} of the virtual bundle
$E-F$. So the Schur functions $S_I(E-F)$ lie in a ring containing the Chern
classes of $E$ and $F$; e.g., we can take the cohomology ring
$H^{\bullet}(X, {\bf Z})$ or the Chow ring $A(X)$ for a smooth,
algebraic $X$. More generally, the Schur functions $S_I(-)$ are well 
defined on the Grothendieck group $K_0(X)$ of vector bundles.
In particular, given a vector bundle $E$ and a partition $I$, we shall
write respectively
$$
S_I(E) \ \ \ \hbox{and} \ \ \ S_I(-E)
$$
for $S_I(E-0)$ and $S_I(0-E)$, where $0$ is the zero vector bundle.

\medskip

We have for any partition $I$\,,
\begin{equation}\label{dual}
S_I(E^*-F^*)=S_{I^{\sim}}(F-E)\,
\end{equation}
where $I^{\sim}$ denotes the {\it dual partition} of $I$. In particular,
$$
S_i(E^*-F^*)=c_i(F-E)
$$
for any $i$, so that the Thom polynomial can be equivalently expressed as
a polynomial in the
$$
S_i(R_m^*-R_n^*)\hbox{'s}\,.
$$
Here, and in the following, we omit pull back indices.
Or, evaluated in the Chern classes $c_i(M)$, $c_j(N)$ of the manifolds 
involved in the map $f: M\to N$, it is expressed as a polynomial in the
$$
S_i(TM^*-f^*TN^*)\hbox{'s}\,.
$$
This convention was used in \cite{P}, \cite{P1} and is used in the present 
paper.

We refer to \cite{L} and \cite{PT} for the theory of Schur functions
(a brief account of Schur functions applied to Thom polynomials,
can be found in \cite[Sect.~3]{P}). We shall use a notation for partitions
indicating the number of times each integer occurs as a part.
For example, we shall denote the partitions $(1,2,2,3)$, $(1,1,1,2,2)$,
and $(1,1,1,1,1,2)$ respectively by $12^23$, $1^32^2$, and $1^52$.

\medskip

Using the theory of supersymmetric functions (cf., e.g., \cite{PT}),
the Thom-Damon theorem can be rephrased by saying that there exist
$\alpha_I\in {\bf Z}$ such that
\begin{equation}\label{alpha}
{\cT}^{\Sigma}=\sum_I \alpha_I S_I(R_m^* - R_n^*)\,,
\end{equation}
the sum is over partitions $I$ with $|I|=\codim(\Sigma)$. The expression in
Eq.~(\ref{alpha}) is unique {\it (loc.cit.)}.

\begin{example}\label{Ep} \ Let $S_I=S_I(R_m^*\moins R_n^*)$. We have the
following three formulas valid for $m-n=0,1,2$, respectively:

\medskip

${\cT}^{A_4}=24S_{4} + 26S_{13} + 10S_{2^2} + 9S_{1^22} + S_{1^4}$,

${\cT}^{III_{2,3}}=8S_{35} + 4S_{134} + 2S_{23^2}$,

${\cT}^{I_{2,2}}=S_{24^2}\plus 3S_{145}\plus 7S_{46}\plus 3S_{5^2}$.

\end{example}
Note that all the coefficients in the formulas in Example \ref{Ep} are
nonnegative. For a more extensive list of examples, also for maps
${\bf C}^2\to {\bf C}$, see the last section.

\section{Cone classes for ample vector bundles}

In the proof of our main result, we shall use the following two results
of Fulton and Lazarsfeld from \cite{FL}. Recall first some classical
definitions and facts from \cite{F} (we shall also follow the notation from
this book). Let $E$ be a vector bundle of rank $e$ on $X$.
By a {\it cone} in $E$ we mean a subvariety of $E$ which is stable under
the natural ${\mathbb G}_{\rm m}$-action on $E$. If $C\subset E$ is a cone 
of
pure dimension $d$, then one may intersect its cycle $[C]$ with the
zero-section of the vector bundle:
\begin{equation}
z(C,E):=s_E^*([C])\in A_{d-e}(X)\,,
\end{equation}
where $s_E^*: A_d(E)\to A_{d-e}(X)$ is the Gysin map determined by the
zero section $X\to E$.
For a projective variety $X$, there is well defined {\it degree}
$\int_X: A_{0}(X)\to {\bf Z}$. We recall first the key technical result
of \cite{FL}.

\begin{theorem}\label{PFL}{\rm (\cite[Theorem 2.1]{FL})}
\label{pr2} Let $E$ be an ample vector bundle of rank $e$ on a
projective variety $X$ of dimension $e$, and let $C\subset E$ be a
cone of pure dimension $e$. Then we have
$$
\int_X \ z(C,E) >0.
$$
\end{theorem}
Under the assumptions of the theorem,
we also have in $H_0(X,{\bf Z})$ the homology analog of $z(C,E)$,
denoted by the same symbol, and the homology degree map
$\deg_X: H_0(X,{\bf Z})\to {\bf Z}$. They are compatible with their
Chow group counterparts via the cycle map:
$A_0(X) \to H_0(X,{\bf Z})$ (cf. \cite[Chap. 19]{F}). Thus we have
\begin{equation}\label{ineq}
\deg_X \bigl(z(C,E)\bigr)>0\,.
\end{equation}

\smallskip

Let $P$ be a symmetric polynomial in $e$ variables and of degree $n$.
It has a unique presentation as a ${\bf Z}$-linear combination
$$
\sum_I {\beta_I}{\cal S}_I\,,
$$
where $|I|=n$, and ${\cal S}_I$ is the classical Schur function
(cf. \cite{L}). Recall that ${\cal S}_I$ is
defined, e.g., by the $e\times e$ determinant (\ref{schur})
with $a$'s replaced by the given variables, and all $b$'s equal to zero.
For any vector bundle $E$ of rank $e$, we define $P(E)$ to be the image of
$P$ under the homorphism which sends the variables to the Chern roots
of $E$. In other words,
$$
P(E)=\sum_I \beta_I S_I(E)\,.
$$
One says that $P$ is {\it numerically positive for ample vector bundles}
if for every projective variety $X$ of dimension $n$, and every ample
vector bundle $E$ of rank $e$ on $X$,
$$
\int_X P(E)
$$
is strictly positive. The second result characterizes polynomials
numerically positive for ample vector bundles, with the help of Schur
functions.

\begin{theorem}{\rm (\cite[Theorem I]{FL})}\label{th2} A homogeneous
polynomial
$$
\sum_I \beta_I {\cal S}_I\,,
$$
where $\beta_I\in {\bf Z}$,
is numerically positive for ample vector bundles iff
for any partition $I$ we have $\beta_I\ge 0$, and additionally
$\sum_I \beta_I>0$.
\end{theorem}
For example, the $n$th Chern class ${\cal S}_{1^n}$ is numerically positive
provided $n\le e$, and in the surface case, i.e. for $n=2$, any polynomial
numerically positive for ample vector bundles is either a positive integer
multiple of ${\cal} S_1$, or a {\bf Z}-linear combination
$$
aS_2+bS_{1^2}\,,
$$
where $a,b\ge 0$ and $a+b>0$ (these results are due respectively to
Bloch-Gieseker and Kleiman, cf. the references in \cite{FL}).

\section{Main result}

We start with some preliminaries that are known in topology
in a much more general framework.
We use the notation from Section \ref{Thom}.
We first pull back the bundle
$\widetilde{\cJ}$ from $BG$ to $BGL_m\times BGL_n$ via the map
induced by the embedding
$$
GL_m \times GL_n \hookrightarrow {\Aut}_m\times {\Aut}_n\,.
$$
Since $GL_m \times GL_n$ acts linearly on $\cJ$, the obtained
pullback bundle is now the  vector bundle on $BGL_m\times BGL_n$
associated with the representation of $Gl_m\times GL_n$ on $\cJ$:
$$
\cJ(R_m,R_n):=\bigl(\oplus_{i=1}^k {\Sym}^i(R_m^*)\bigr) \otimes R_n\,.
$$
The bundle $\cJ(R_m,R_n)$ contains the preimage of
$\widetilde{\Sigma}\subset \widetilde{\cJ}$, denoted by $\Sigma(R_m,R_n)$,
whose dual class is given by the RHS of Eq.~(\ref{alpha}).

Consider, more generally, a pair of vector bundles $E$ and $F$ of
ranks $m$ and $n$ on a variety $X$. We define the following
vector bundle on $X$:
$$
\cJ(E,F):=\bigl(\oplus_{i=1}^k {\Sym}^i(E^*)\bigr) \otimes F\,.
$$
In fact, the pair of bundles $(E,F)$ corresponds to a principal
$GL_m\times GL_n$-bundle $B(E,F)$ and
$$
\cJ(E,F)=B(E,F)\times_{GL_m\times GL_n}\cJ
$$
is the bundle associated with the representation. Similarly, we define
the singularity set
$$
\Sigma(E,F):=B(E,F)\times_{GL_m\times GL_n}\Sigma\subset\cJ(E,F)\,.
$$
The dual class\footnote{A precise definition of the dual class of
$[\Sigma(E,F)]$, in the case when $X$ is singular, is given in Note
\ref{note}.} of $[\Sigma(E,F)]$ in
$$
H^{2\codim(\Sigma)}(\cJ(E,F), {\bf Z})\cong H^{2\codim(\Sigma)}(X, {\bf Z})
$$
is equal to
\begin{equation}\label{alpha''}
\sum_I \alpha_I S_I(E^*-F^*)\,,
\end{equation}
where the $\alpha_I$'s were defined in Eq.~(\ref{alpha}). The
argument for that is fairly standard but one has to pass to
topological homotopy theory, where each pair of bundles can be
pulled back from the universal pair $(R_m, R_n)$ of bundles on
$BGL_m\times BGL_n$. It is possible to work entirely with the
algebraic varieties. One can use the Totaro construction and
representability for affine varieties (\cite[proof of Theorem
1.3]{To}).

The main result of the present paper, suggested/conjectured in \cite{P}, 
\cite {P1}, and in \cite {FK} for Thom-Boardman singularities, is

\begin{theorem}\label{Tpos} Let $\Sigma$ be a stable, nontrivial class of
singularities. Then the Thom polynomial ${\cT}^{\Sigma}$ is nonzero,
and for any partition $I$ the coefficient $\alpha_I$ in the Schur
function expansion of the Thom polynomial
$\cT^{\Sigma}$ (cf. Eq.~(\ref{alpha})) is nonnegative.
\end{theorem}
\proof
We follow the notation from the first part of this section.
Let $c=\codim(\Sigma)$ (so that for any partition $I$ with $\alpha_I\ne 0$,
we have $|I|=c$). Suppose that $X$ is a projective variety of dimension $c$,
$F$ is an ample vector on $X$ with $\rank(F)=n'=n+r \ge c$, and
$E={\bf 1}_X^{m'}$ is a trivial bundle of rank $m'=m+r$ (for some $r\ge 0$).
The variety $\Sigma({\bf 1}^{m'},F)$ is a cone in $\cJ({\bf 1}^{m'},F)$
because ${\mathbb G}_{\rm m} \subset \Aut_{n'}$.
From the first part of this section, we know that the cone class
$$
z\bigl(\Sigma(E,F),{\cJ}(E,F)\bigr)\in H_0(X,{\bf Z})
$$
is dual to the universal expression (\ref{alpha''}):
$$
\sum_I \alpha_I S_I(E^*-F^*)=\sum_I \alpha_I S_I(-F^*)
=\sum_I \alpha_I S_{I^{\sim}}(F)
$$
(we use here Eqs. (\ref{schur}) and (\ref{dual})).
Using the notation of Theorem \ref{th2}, we set
$$
P=\sum_I \alpha_I {\cal S}_{I^{\sim}}\,.
$$
Since a direct sum of ample vector bundles is ample
\cite[Proposition 2.2]{H}, the vector bundle
$$
\cJ({\bf 1}^{m'},F)=F^{\oplus N}
$$
(for some integer $N$) is ample. Hence we have we have by the
inequality (\ref{ineq})
$$
\int_X P(F) = \deg_X \bigl(z(\Sigma({\bf 1}^{m'},F), F^{\oplus N})\bigr)>0\,,
$$
that is, the polynomial $P$ is numerically positive for ample vector bundles.
We conclude, by Theorem \ref{th2}, that the coefficients $\alpha_I$
in Eq. (\ref{alpha}) are nonnegative. Moreover, we have
$\sum_I \alpha_I>0$, which implies that ${\cal T}^{\Sigma}\ne 0$.

\medskip

The theorem has been proved.
\qed

\bigskip

Theorems \ref{PFL} and Theorem \ref{th2} are valid over any algebraically
closed field of arbitrary characteristic (cf. \cite{FL}). Therefore our main
result would hold in that generality provided that one develops a
suitable theory of Thom polynomials. Of course, usual cohomology should be
replaced by Chow rings, cf. \cite{To}. By Theorem 1.3 {\it (loc.cit.)},
the theory of characteristic classes is reduced to the calculus in the
Chow ring of a classifying space.

\begin{remark}\rm\label{uwaga}
The functor  of $k$-jets
$$
E,F\mapsto {\cJ}(E,F)=\left(\oplus_{i=1}^k {\Sym}^iE^*\right)\otimes F
$$
and the cone $\Sigma(E,F)$ can be replaced by more general construction:
an arbitrary functor $\phi(E,F)$ and a cone bundle
$\Sigma(E,F)\subset\phi(E,F)$. We assume the following two properties:
\begin{enumerate}
\item the class $z(\Sigma(E,F),\phi(E,F))$ is stable under simultaneous 
addition of the same bundle;
\item the assignment $F\mapsto \phi({\bf 1}^m,F)$ (or $E\mapsto
\phi(E^*,{\bf 1}^n)$) preserves ampleness.
\end{enumerate}
Examples of functors preserving ampleness for fields of characteristic
zero are {\it polynomial functors}. They are, at the same time, quotient
functors and subfunctors of the tensor power functors (cf. \cite{H}).
Theorem \ref{Tpos} remains then valid.
\end{remark}

\begin{note}\label{note} \rm
We give here a precise definition of the dual class of $[\Sigma(E,F)]$
for possibly singular $X$. This class will lie in $H^{2c}(X,{\bf Z})$
(recall that $c=\codim(\Sigma)$).
Let $Y\subset\Sigma$ be the set of singular points of $\Sigma$, i.e.
$$
\Sigma\setminus Y\subset \cJ\setminus Y
$$
is a submanifold. The set $Y$ is $GL_m\times GL_n$-invariant.
For a pair of bundles we construct the fibering
$$
Y(E,F)=B(E,F)\times_{GL_m\times GL_n} Y\subset \Sigma(E,F)\subset \cJ(E,F)
$$
with fiber $Y$. We have $\codim_{\bf C}Y \ge c+1$ and therefore
$$
H^i(\cJ\setminus Y,{\bf Z})\cong H^i(\cJ,{\bf Z})
$$
for $i\leq 2c$. Hence
$$
H^i(\cJ(E,F)\setminus Y(E,F), {\bf Z})\cong H^i(\cJ(E,F),{\bf Z})
$$
in the same range of degrees, and it is enough to define the desired
class in $H^{2c}(\cJ(E,F)\setminus Y(E,F),{\bf Z})$. Note that
$\Sigma(E,F)\setminus Y(E,F)$ has a normal bundle. Its Thom class
defines an element in
$$
H^{2c}(\cJ(E,F)\setminus Y(E,F),\cJ(E,F)\setminus\Sigma(E,F);{\bf Z})\,.
$$
The image of this element in
$$
H^{2c}(\cJ(E,F)\setminus Y(E,F), {\bf Z})\cong H^{2c}(\cJ(E,F),{\bf Z})
\cong H^{2c}(X,{\bf Z})
$$
is the desired class. Verification that this class is natural with respect
to pull back of vector bundles is left to the reader.
\end{note}

\section{Examples}

Let $S_I=S_I(R_m^*\moins R_n^*)$ in the notation of Section 2.
We list several examples of singularities $\Sigma[n\moins m]: M \to N$,
where $m=\dim(M) \le n=\dim(N)$. All of them were computed by the ``method
of restriction equations'' of \cite{Rim}.
In the first paragraph we give the Schur function expansions of Thom
polynomials for singularities $\Sigma[0]$ with codimension $\le 6$
from \cite{Rim}, p.~508. In the second paragraph, we give the Schur
function expansion of the Thom polynomials for the singularities
$\Sigma[1])$ from \cite{Rim}, p.~512. In the third paragraph, we list
some examples from \cite{P}, \cite{P1}. (In \cite{P}, \cite{P1} we used 
a ``shifted'' notation: $\Sigma(n-m+1)$.)

\medskip

\parindent 0pt

$A_1[0]$: $S_{1}$

$A_2[0]$: $2S_{2} + S_{1^2}$

$A_3[0]$: $6S_{3} + 5 S_{12} + S_{1^3}$

$A_4[0]$: $24S_{4} + 26S_{13} + 10S_{2^2} + 9S_{1^22} + S_{1^4}$

$I_{2,2}[0]$: $S_{2^2}$

$A_5[0]$: $120S_{5} + 154S_{14} + 92S_{23} + 71S_{1^23} + 14S_{1^32} +
35S_{12^2} + S_{1^5}$

$I_{2,3}[0]$: $4S_{23} + 2S_{12^2}$

$A_6[0]$: $720S_6 + 1044S_{15} + 770S_{24} + 266S_{3^2} + 580S_{1^24} +
455S_{123} + 70S_{2^3} + 155S_{1^33} + 84S_{1^22^2} + 20S_{1^42} +
S_{1^6}$

\smallskip

$I_{2,4}[0]$: $16S_{24} + 4S_{3^2} + 12S_{123} + 5S_{2^3} + 2S_{1^22^2}$

$I_{3,3}[0]$: $2S_{24} + 6S_{3^2} + 3S_{123} + S_{1^22^2}$

\bigskip

$A_1[1]$: $S_2$

$A_2[1]$: $4S_4 + 2S_{13} + S_{2^2}$

$A_3[1]$: $36S_6 + 30S_{15} + 19S_{24} + 5S_{3^2} + 6S_{1^24} + 5S_{123}
+ S_{2^3}$

$A_4[1]$: $507S_8 + 555S_{17} + 391S_{26} + 240S_{35} + 76S_{4^2} +
216S_{1^26} + 210S_{125} + 104S_{134} + 55S_{2^24} + 21S_{23^2} + 24S_{1^35} +
26S_{1^224} + 10S_{1^23^2} + 9S_{12^23} + S_{2^4}$

\smallskip

$III_{2,2}[1]$: $S_{3^2}$

$I_{2,2}[1]$: $3S_{34} + S_{13^2}$

$III_{2,3}[1]$: $8S_{35} + 4S_{134} + 2S_{23^2}$

\bigskip

$I_{2,2}[1]$: $S_{13^2}+3S_{34}$

$I_{2,2}[2]$: $S_{24^2}\plus 3S_{145}\plus 7S_{46}\plus 3S_{5^2}$

$I_{2,2}[3]$: $S_{35^2}\plus 3S_{256}\plus 7S_{157}\plus 3S_{16^2}\plus 15S_{58}
\plus 10S_{67}$

\bigskip\medskip

\noindent
In this paragraph, we list Thom polynomials
for the functions ${\bf C}^2\to {\bf C}$ of the singularities
$A, D, E$ with codimension at most 6, computed by the ``method of restriction
equations'' of \cite{Rim}. The results agree with those in \cite{Ka} given
in another basis.

\medskip

$A_1$: $S_{1^2}$

$A_2$: $2S_{1^3}+2S_{12}$

$A_3$: $5S_{1^4}+11S_{1^22}+6S_{2^2}+6S_{13}$

$A_4$: $12S_{1^5}+44S_{1^32}+44S_{12^2}+56S_{1^23}+36S_{23}+24S_{14}$

$D_4$: $S_{1^5}+3S_{1^32}+6S_{12^2}+2S_{1^23}+4S_{23}$

$A_5$: $30S_{1^6}+160S_{1^42}+248S_{1^22^2}+338S_{1^33}+434S_{123}+328S_{1^24}
+108S_{3^2}+228S_{24}+120S_{15}$

\smallskip

$D_5$: $4S_{1^6}+18S_{1^42}+42S_{1^22^2}+26S_{1^33}+64S_{123}+12S_{1^24}
+24S_{3^2}+24S_{24}$

\smallskip

$A_6$: $79S_{1^7}+566S_{1^52}+1238S_{1^32^2}+1723S_{1^43}+
3473S_{1^223}+2736S_{1^34}+1834S_{13^2}+3898S_{124}+2220S_{1^25}+
1260S_{34}+1632S_{25}+720S_{16}$

\smallskip

$D_6$: $8S_{1^7}+50S_{1^52}+138S_{1^32^2}+118S_{1^43}+
348S_{1^223}+124S_{1^34}+224S_{13^2}+320S_{124}+48S_{1^25}+
144S_{34}+96S_{25}$

\smallskip

$E_6$: $3S_{1^7}+18S_{1^52}+54S_{1^32^2}+39S_{1^43}+
129S_{1^223}+36S_{1^34}+102S_{13^2}+102S_{124}+12S_{1^25}+
60S_{34}+24S_{25}$

\bigskip

\noindent
{\bf Acknowledgments} \ The first author thanks Alain Lascoux for helpful
discussions on Schur functions and Thom polynomials.
The second author thanks Maxim Kazarian for suggesting Remark \ref{uwaga}.
Both authors thank Ozer Ozturk for his help with computer computations.

\end{document}